\documentclass[11pt,amsfonts]{article}
\usepackage{graphicx}
\usepackage{latexsym}
\usepackage{amssymb}
\usepackage{amsmath}
\usepackage{layout}
\newtheorem{prop}{Proposition}

\newtheorem{theorem}{Theorem}
\newtheorem{remark}{Remark}

\def\real{{\mathord{{\rm I\kern-2.8pt R}}}}        
\def\inte{{\mathord{{\rm I\kern-2.8pt N}}}}

\def\sZZ{{\rm Z\kern-2.8ptem{}Z}}

\def\z{{\mathchoice
  {\sZZ}
  {\sZZ}
  {\rm Z\kern-0.30em{}Z}
  {\rm Z\kern-0.25em{}Z} }}
\def\sQQ{{\kern 0.27em \vrule height1.45ex width0.03em depth0em
          \kern-0.30em \rm Q}}
\def\qu{{\mathchoice
    {\sQQ}
    {\sQQ}
  {\kern 0.225em \vrule height1.05ex width0.025em depth0em \kern-0.25em \rm Q}
  {\kern 0.180em \vrule height0.78ex width0.020em depth0em \kern-0.20em \rm Q}
        }}
\def\sCC{{\kern 0.27em \vrule height1.45ex width0.03em depth0em
          \kern-0.30em \rm C}}
\def\complex{{\mathchoice
    {\sCC}
    {\sCC}
  {\kern 0.225em \vrule height1.05ex width0.025em depth0em \kern-0.25em \rm C}
  {\kern 0.180em \vrule height0.78ex width0.020em depth0em \kern-0.20em \rm C}
        }}

\newcommand{\ba}{\begin{array}}
\newcommand{\ea}{\end{array}}
\newcommand{\be}{\begin{equation}}
\newcommand{\ee}{\end{equation}}
\newcommand{\bea}{\begin{eqnarray}}
\newcommand{\eea}{\end{eqnarray}}
\newcommand{\beaa}{\begin{eqnarray*}}
\newcommand{\eeaa}{\end{eqnarray*}}

\def\z{\zeta}

\font\tenmath=msbm10 \font\sevenmath=msbm7 \font\fivemath=msbm5
\newfam\mathfam \textfont\mathfam=\tenmath
\scriptfont\mathfam=\sevenmath \scriptscriptfont\mathfam=\fivemath

\def \={{\buildrel {\rm (law)} \over =}}

\def\qed{ \hfill \vrule width.25cm height.25cm depth0cm\smallskip}

\newcommand{\basa}{\begin{assumption}}
\newcommand{\easa}{\end{assumption}}

\newcommand{\bas}{\begin{assum}}
\newcommand{\eas}{\end{assum}}

\newcommand{\ignore}[1]{}
\begin{document}

\renewcommand{\thefootnote}{\fnsymbol{footnote}}

\title{Donsker type theorem for the Rosenblatt process and a binary market model }

\vskip1cm

\author{Soledad Torres  $^{1}\quad$Ciprian A. Tudor $^{2}\vspace*{0.1in}$\\$^{1}$
 Departamento de Estad\'istica, Universidad de Valpara\'iso,\\ Casilla
 123-V,
4059 Valparaiso, Chile.
\\soledad.torres@uv.cl\vspace*{0.1in}\\$^{2}$SAMOS/MATISSE,
Centre d'Economie de La Sorbonne,\\ Universit\'e de
Panth\'eon-Sorbonne Paris 1,\\90, rue de Tolbiac, 75634 Paris Cedex
13, France.\\tudor@univ-paris1.fr\vspace*{0.1in}}  \maketitle

 \begin{abstract}
In this paper, we prove a Donsker type approximation theorem for the
Rosenblatt process, which is a selfsimilar stochastic process
exhibiting long range dependence. By using numerical results and
simulated data, we show that this approximation performs very well.
We use this result to construct a binary market model driven by this
process and we show that the model admits arbitrage opportunities.
\end{abstract}

 \vskip0.5cm

{\bf  2000 AMS Classification Numbers: }60F17, 91B70.

 \vskip0.3cm


\vskip0.3cm

\section{Introduction}
Long range dependence  stochastic processes have been intensively
used  as models for different physical phenomena. First, these
properties appeared in empirical studies  in areas like hydrology
and geophysics; more recently, they appeared to play an important
role in network traffic analysis and telecommunications. As a
consequence,  efficient mathematical models based on long range
dependence (or long memory) processes have been proposed in these
directions.

The notions of long range dependence and  selfsimilarity have also
been considered in mathematical finance. An excellent survey on the
different aspects of the appearance of the long range dependence in
practice is the paper \cite{rama}. The debate on the presence of
long memory in stock prices is actually not new. The  idea that
asset returns could exhibit long range dependence comes from
Mandelbrot (\cite{Mand}) and then observed in several empirical
studies. We refer, among others,  to \cite{WTT} for concrete
examples and  for interesting comments on this question. We also
mention that some other authors rejected the idea of the presence of
long memory in asset returns (see e.g. \cite{Lo}).

A rather general opinion is that long range dependence in financial
models is strongly related to the presence of arbitrage. For
example, in the case of market  models  driven by the fractional
Brownian motion, this has been explicitly shown by Rogers
(\cite{Rog}) or Sottinen (\cite{Sot}). In special situations, for
example under transaction costs, arbitrage could be eliminated (see
\cite{Gua} or \cite{Salo}). Different approaches, based on
Wick-It\^o calculus, have been developed in e.g. \cite{Ok1},
\cite{Ok2}. Models driven by long range dependence processes others
than the fractional Brownian motion have been, from the stochastic
calculus point of view, less considered. This is actually one of the
motivations of our work: we propose a binary market model driven by
a non-Gaussian selfsimilar
 process (called the Rosenblatt
 process) which exhibits long range dependence  and we show that
 this type of model admits arbitrage opportunities. It is actually
 known that in general arbitrage-free  models imply the fact that the
 price process is a semimartingale (\cite{DS}) but it has been recently proved
 in \cite{CoRu} that by modifying the class of admissible
 strategies, one can also consider  arbitrage-free  models driven by
 non-martingales.

 To construct our binary market model based on the Rosenblatt
 process  we need a Donsker type theorem to approximate in law this
 process by some disturbed two-dimensional random  walks; this
 results could be useful by itself. In fact, this theorem extends   a result by
 Sottinen \cite{Sot} and represents a variant of the so-called Non
 Central Limit Theorem proved in \cite{DM} and \cite{Taqqu}. We
 mention that, since we are now in a non-Gaussian context, the
 proof of this result demands different techniques.

 Our paper is organized as follows. In Section 2 we describe the
 basic properties of the Rosenblatt process. Section 3 contains the
 proof of the Donsker theorem to approximate weakly, in the
 Skorohod topology, the Rosenblatt  process by walks. In Section 4,
 we introduce our binary  market model which is showed to converge
 to the Black and Scholes model with Rosenblatt noise. We show that
 the model admits arbitrage opportunities and we construct a such
 opportunity.

Finally, in Appendix  we present numerical results based on
simulated data which show that the approximation method performs
very well. Related numerical results can be found in \cite{PA} or
\cite{Pi}.

 \section{Preliminaries}
 The Rosenblatt process appears as a limit in the so called Non
 Central Limit Theorem (see \cite{DM} or \cite{Taqqu}). We recall
 the general context. Consider $(\xi _{n}) _{n\in \mathbb{Z}}$ a stationary Gaussian
sequence with mean zero and covariance $1$ such that its correlation
function satisfies
\begin{equation}
\label{corr} r(n) := \mathbf{E}\left( \xi _{0} \xi _{n} \right) =
n^{\frac{2H-2} {k}} L(n),
\end{equation}
with $k\geq 1$ integer, $H\in (\frac{1}{2}, 1) $  and $L$ is a
slowly varying function at infinity . Denote by $H_{m}(x)$ the
Hermite polynomial of degree $m$ given by $ H_{m}(x)=
(-1)^{m}e^{\frac{x^{2}}{2}}\frac{d^{m}}{dx^{m}}e^{-\frac{x^{2}}{2}}
$ .  Let $g $ be a function such that $\mathbf{E}(g(\xi _{0} ))=0$
and $\mathbf{E}(g(\xi _{0})^{2})<\infty. $ Suppose that $g$ {\em has
Hermite rank } equal to $k$; that is, if $g$ admits the following
expansion in Hermite polynomials
\begin{equation*}
g(x)= \sum _{j\geq 0} c_{j} H_{j}(x), \hskip0.5cm c_{j}=
\frac{1}{j!} \mathbf{E}\left( g(\xi _{0} H_{j} (\xi _{0}))
\right),
\end{equation*}
then
\begin{equation*}
k=\min \{ j; c_{j}\not= 0\}.
\end{equation*}
 Since $\mathbf{E}\left[ g(\xi
 _{0})\right] =0$,  we have $k\geq 1$.  {\em The Non Central Limit
 Theorem} (\cite{DM}, \cite{Taqqu}) says that the sequence of
 stochastic processes
 \begin{equation*}
 \frac{1}{n^{H}} \sum _{ j=1} ^{[nt]} g(\xi_{j})
 \end{equation*}
 converges as $n \to \infty $, in the sense of finite dimensional
 distributions,  to the process (called {\em the Hermite process})
 \begin{equation}
\label{hermite} Z^{k}_{H} (t)= c(H,k)
\int_{\mathbb{R}^{k}}\int_{0}^{t} \left( \prod _{j=1}^{k} (s-y_{i}
)_{+} ^{ -\left( \frac{1}{2} + \frac{1-H}{k} \right) }  \right) ds
dB(y_{1} )\ldots dB(y_{k}),
\end{equation}
where $x_{+}=\max(x,0)$ and the above integral is a multiple
Wiener-It\^o stochastic integral with respect to a Brownian motion
$B(y))_{y\in \mathbb{R}}$ (see \cite{N}) .

Let us list some basic properties of the Hermite processes.
\begin{description}
\item{$\bullet$ } it exhibits long-range
dependence (the covariance function decay at a power function at
zero -"Joseph effect")

\item{$\bullet$ }it is $H$-selfsimilar in the sense that for any $c>0$,
$(Z^{k}_{H}(ct)) =  ^{(d)} (c^{H} Z^{k}_{H} (t))$, where $ " =^{(d)}
"$ means equivalence of all finite dimensional distributions

\item{$\bullet$ }it
has stationary increments, that is, the joint distribution of
$(Z^{k}_{H}(t+h)-Z^{k}_{H}(h), t\in [0,T])$ is independent of $h>
0$.

\item{$\bullet$ } the covariance function is
$$\mathbf{E}(Z^{k}_{H}(t)Z^{k}_{H}(s)) = \frac{1}{2}\left( t^{2H} + s^{2H } -\vert
t-s \vert ^{2H} \right), \hskip0.5cm s,t\in [0,T] $$

and consequently, for every $s,t\in [0,T]$
\begin{equation}
\label{zt-zs} \mathbb{E}\left| Z_{H}^{k}(t) -Z_{H}^{k}(s)\right|
^{2} = \vert t-s\vert ^{2H}
\end{equation}

\item{$\bullet$ }  the Hermite process is Hold\"er continuous of order $\delta <H$

\item{$\bullet$ } if $k\geq 2$, then $Z^{k}_{H}$ is non-Gaussian.

\end{description}

When $k=1$ the process given by (\ref{hermite}) is nothing else that
the {\em fractional Brownian motion (fBm) } with Hurst parameter
$H\in (0, \frac{1}{2})$. For $k\geq 2$ the process is not Gaussian.
If $k=2$ then the process (\ref{hermite}) is known as {\em the
Rosenblatt process }(it has been actually named by M. Taqqu).

We focus here  our attention on the case $k=2$.  We will work with
the representation of this processes as integral with respect to a
Wiener process on a finite interval. Recall
 that the fBm $B^{H}$ with Hurst parameter  $H>\frac{1}{2}$ can be written as
 \begin{equation}
 \label{B1}
B_{t}^{H}= \int_{0}^{t} K^{H}(t,s)dW_{s}, \hskip0.5cm t\in [0,T]
 \end{equation}
 with $(W_{t}, t\in [0,T])$ a standard Wiener process and
 \begin{equation}
 \label{K}
K^{H}(t,s)= c_{H} s^{\frac{1}{2}-H} \int _{s}^{t}
(u-s)^{H-\frac{3}{2}} u^{H-\frac{1}{2}}  du
 \end{equation}
 where $t>s$ and $c_{H} =\left( \frac{ H(2H-1) }{\beta( 2-2H, H-\frac{1}{2}) } \right)
^{\frac{1}{2}}.$  From (\ref{K}) we obtain that for $t>s$,
\begin{equation}
\label{dK} \frac{\partial K}{\partial t} (t,s)= c_{H}\left(
\frac{s}{t} \right) ^{\frac{1}{2}-H} (t-s)^{H-\frac{3}{2}}.
\end{equation}
Aa analogous representation for the Rosenblatt process has been
given in \cite{CAT}. We have
\begin{equation}
\label{rose2} Z^{2}_{H}(t):=Z_{t}=^{(d)} d(H) \int
_{0}^{t}\int_{0}^{t} \left[ \int_{ y_{1} \vee y_{2} }^{t}
\frac{\partial K^{H'} }{\partial u} (u,y_{1} )  \frac{\partial
K^{H'} }{\partial u} (u,y_{2} )du \right]dW_ {y_{1}} dW_{y_{2}}
\end{equation}
where $(W_{t}, t \in [0,T])$ is a Brownian motion,
\begin{equation}
\label{H'} H'=\frac{H+1}{2}
\end{equation}
and $ d(H)= \frac{1}{H+1} \left( \frac{H}{2(2H-1)} \right)
^{-\frac{1}{2}}.$  Note that $H>\frac{1}{2}$ implies
$H'>\frac{3}{4}$.

\section{Convergence in law to the Rosenblatt process}

This part in consecrated to a Donsker invariance principle for the
Rosenblatt process. From now on, we will consider the Rosenblatt
process to be given by the formula (\ref{rose2}). We will denote,
for every $t\in [0,T]$
\begin{equation*}
F(t,y_{1}, y_{2}) =d(H)1_{[0,t]}(y_{1})1_{[0,t]}(y_{2})\int_{ y_{1}
\vee y_{2} }^{t} \frac{\partial K^{H'} }{\partial u} (u,y_{1} )
\frac{\partial K^{H'} }{\partial u} (u,y_{2} )du
\end{equation*}
and then
\begin{equation}
\label{rose3} Z_{t}= \int _{0} ^{T} \int _{0} ^{T} F(t, y_{1},
y_{2}) dW(y_{1}) dW(y_{2}), \hskip0.5cm t\in [0,T].
\end{equation}
The kernel $K^{H'}$ (denoted simply by $K$ in the sequel) is  the
standard kernel of the fractional Brownian motion (\ref{K}).

\vskip0.3cm

 Let us first recall some known facts. Consider $(\xi_i) _{i\geq 1}$ a sequence of i.i.d random
variables with $\mathbb{E}(\xi _i)= 0$ and $\mathbb{E}(\xi_i^2) =1$.
The Donsker Invariance Principle says that the sequence of processes
\begin{equation}
\label{rw} W^{n}_{t}= \frac{1}{\sqrt{n}} \sum _{i=1}^{[nt]} \xi _{i}
\end{equation}
converges weakly, in the Skorohod topology, to a standard Brownian
motion. Here $[x]$ denotes the biggest integer smaller than $x$.

This result has been extended in \cite{Sot} to the fractional
Brownian motion (see also \cite{Nie}). Define
$$K ^{n} (t,s) :=n\int_{s-\frac{1}{n}}^{s} K(\frac{[nt]}{n},u)du, \hskip0.5cm n\geq 1 $$
and put
\begin{equation*}
B^{n}_{t}= \int_{0}^{t}K^{n}(t,s) dW^{n}_{s} = \sum _{i=1}^{[nt]}
n\int _{\frac{i-1}{n}}^{\frac{i}{n}}K(\frac{[nt]}{n}, s)ds \frac{\xi
_{i} }{\sqrt{n}}, \hskip0.5cm n\geq 1 .
\end{equation*}
Then it has been proved in \cite{Sot} that the disturbed random walk
$B^{n}$ converges weakly to the fractional Brownian motion.

From the above results and the representation (\ref{rose3}) it is
quite natural to  define the following approximation for the
Rosenblatt process
\begin{equation}
\label{zn} Z^{n}_{t}= \sum_{i,j=1; i\not= j}^{[nt]}n^{2}
\int_{\frac{i-1}{n} }^{\frac{i}{n}}\int_{\frac{j-1}{n}}^
{\frac{j}{n}}F\left( \frac{[nt]}{n}, u,v\right) dvdu \frac{\xi
_{i}}{\sqrt{n}}\frac{\xi _{j}}{\sqrt{n}}, \hskip0.5cm t\in [0,T].
\end{equation}

\begin{remark}
We eliminate the diagonal "i=j" because the Rosenblatt process is
defined as a double Wiener-It\^o integral and as a consequence it
has zero mean. When the diagonal $i=j$ is included  in the sum
(\ref{zn}) then the limit is in general a double Stratonovich
integral (see \cite{HM} or \cite{SU}).

\end{remark}

\vskip0.5cm

\begin{prop}\label{findim}
 The family of stochastic processes $(Z^{n}_{t})_{t\in [0,T]} $
 converges in  the sense of finite dimensional distributions to the
 process $(Z_{t})_{t\in [0, T]}$ (\ref{rose3}).
\end{prop}
{\bf Proof: } We will proof this result in several steps.

\vskip0.3cm

{\it Step 1: } Let us consider an arbitrary sequence of partitions
of the interval $[0,T]$ of the form
$$ \pi ^{m}: 0=t_{0}^{m}< t_{1}^{m} <\ldots < t_{m}^{m}=T$$
with $\vert \pi ^{m} \vert \to 0$ as $m\to \infty$. Define
\begin{equation}
\label{zpim} Z^{\pi ^{m}}_{t}= \sum _{i,j=1; i\not= j} ^{m}
\frac{1}{ \vert \Delta _{i}^{m} \vert \vert \Delta _{j}^{m} \vert }
\left( \int _{\Delta _{i}^{m}} \int _{\Delta _{j}^{m}}F(t,u,v) dvdu
\right) W(\Delta _{i} ^{m}) W(\Delta _{j} ^{m})
\end{equation}
where we denoted by $\Delta _{i} ^{m} = [t_{i-1}^{m}, t_{i}^{m}) $
and by
$$W (\Delta _{i}^{m})= W_{t_{i}^{m}} -W_{t_{i-1}^{m}}.$$
Then it follows from \cite{SU}, Theorem 3.4, or \cite{HM} that for
fixed $t$ the sequence $Z_{t}^{\pi ^{m}} $ converges in $L^{2}
(\Omega )$ as $\vert \pi ^{m}\vert \to 0$ to the multiple
Wiener-It\^o integral of $F(t, \cdot )$ with respect to the Brownian
motion $W$
 $$\int_{0}^{T}\int_{0}^{T}F(t,u,v) dW_{u} dW_{v} =Z_{t}.$$

\vskip0.3cm

{\it Step 2: } Secondly, define the process
\begin{equation}
\label{zpimn} Z^{\pi ^{m}, n} _{t}=\sum _{i,j=1; i\not= j} ^{m}
\frac{1}{ \vert \Delta _{i}^{m} \vert \vert \Delta _{j}^{m} \vert }
\left( \int _{\Delta _{i}^{m}} \int _{\Delta _{j}^{m}}F(t,u,v) dvdu
\right) W^{n}(\Delta _{i} ^{m}) W^{n}(\Delta _{j} ^{m})
\end{equation}
where $W^{n}$ is the random walk given by (\ref{rw}). Then clearly,
for fixed $m$, as $n$ goes to $\infty$, the finite dimensional
distributions of $Z^{\pi ^{m}, n} $ converges to the finite
dimensional distributions of $Z^{\pi ^{m}}$ (this comes from the
weak convergence of $W^{n}$ to the Wiener process $W$).

\vskip0.3cm

{\it Step 3: } We prove now that for every $t\in [0,T]$, the
sequence  $Z^{\pi ^{m}, n} _{t}$ converges in $L ^{2} (\Omega )$ to
$Z^{', n}_{t}$ as $m\to \infty$, where
\begin{equation}
\label{zprime} Z^{', n}_{t}=\sum_{i,j=1; i\not= j}^{[nt]}n^{2}
\int_{\frac{i-1}{n} }^{\frac{i}{n}}\int_{\frac{j-1}{n}}^
{\frac{j}{n}}F\left(t, u,v\right) dvdu \frac{\xi
_{i}}{\sqrt{n}}\frac{\xi _{j}}{\sqrt{n}}, \hskip0.5cm t\in [0,T].
\end{equation}
Consider the sequence
 $$F^{\pi ^{m}}(t, u,v) = \sum _{i,j=1; i\not= j} ^{m}
\frac{1}{ \vert \Delta _{i}^{m} \vert \vert \Delta _{j}^{m} \vert }
\left( \int _{\Delta _{i}^{m}} \int _{\Delta _{j}^{m}}F(t,u,v) dvdu
\right) 1_{\Delta _{i}^{m}}(u)  1_{\Delta _{j}^{m}}(v).$$
 Then $F^{\pi ^{m}}(t, \cdot )$ converges to $ F(t, \cdot )$ in $L^{2} ([0,T]^{2})$
as $m \to \infty$ (see \cite{SU}, \cite{HM}).

First note that $Z^{', n}_{t}$ can be approximated in
$L^{2}(\Omega)$ as $m \to \infty $ by
$$ Z^{',\pi ^{m}, n} _{t}=  \sum _{k,l=1; k\not=l}^{[nt]} n^{2}\int _{
\frac{k-1}{n}} ^{\frac{k}{n}} \int  _{
\frac{l-1}{n}}^{\frac{l}{n}}F^{\pi ^{m}}(t,u,v) dvdu \frac{\xi
_{k}}{\sqrt{n}} \frac{\xi_{l}}{\sqrt{n}}.$$ Indeed,
\begin{eqnarray*}
&&\mathbb{E}\left| Z^{',\pi ^{m}, n} _{t} - Z^{', n}_{t} \right| ^{2}\\
&=&\sum _{k,l=1; k\not=l}^{[nt]}n^{2} \left(     \int _{
\frac{k-1}{n}} ^{\frac{k}{n}} \int  _{
\frac{l-1}{n}}^{\frac{l}{n}}\left( F^{\pi
^{m}}(t,u,v)-F(t,u,v)\right)  dvdu \right)
^{2}\mathbb{E}(\xi _{k}^{2}) \mathbb{E}(\xi_{l}^{2})\\
&\leq &\sum _{k,l=1; k\not=l}^{[nt]} \int _{ \frac{k-1}{n}}
^{\frac{k}{n}} \int  _{ \frac{l-1}{n}}^{\frac{l}{n}}(F^{\pi
^{m}}(t,u,v)-F(t,u,v) ) ^{2}dvdu\\
&\leq& \int_{0}^{T} \int_{0}^{T}(F^{\pi ^{m}}(t,u,v)-F(t,u,v) )
^{2}dvdu
\end{eqnarray*}
and  this clearly goes to zero as $m\to \infty$.

It remains to observe that $Z^{', \pi ^{m}, n} _{t}$ is equal  to
$Z^{\pi ^{m}}_{t}$ for every $t,m,n$. We can write, if $\lambda $
denotes the Lebesque measure,
\begin{eqnarray*}
Z^{', \pi ^{m}, n}_{t} &=& \sum _{k,l=1; k\not=l}^{[nt]} n^{2}\int
_{ \frac{k-1}{n}} ^{\frac{k}{n}} \int  _{
\frac{l-1}{n}}^{\frac{l}{n}} dvdu\\
&&\times \left( \sum _{i,j=1; i\not= j} ^{m} \frac{1}{ \vert \Delta
_{i}^{m} \vert \vert \Delta _{j}^{m} \vert } \left( \int _{\Delta
_{i}^{m}} \int _{\Delta _{j}^{m}}F(t,x,y) dydx \right) 1_{\Delta
_{i}^{m}}(u) 1_{\Delta _{j}^{m}}(v) \right)  \frac{\xi
_{k}}{\sqrt{n}} \frac{\xi_{l}}{\sqrt{n}}\\
&=&\sum _{i,j=1; i\not= j} ^{m} \frac{1}{ \vert \Delta _{i}^{m}
\vert \vert \Delta _{j}^{m} \vert } \left( \int _{\Delta _{i}^{m}}
\int _{\Delta _{j}^{m}}F(t,x,y) dydx \right) \\
&&\times \sum _{k,l=1; k\not=l}^{[nt]} n^{2}\frac{\xi
_{k}}{\sqrt{n}} \frac{\xi_{l}}{\sqrt{n}}\int _{ \frac{k-1}{n}}
^{\frac{k}{n}} \int _{ \frac{l-1}{n}}^{\frac{l}{n}} 1_{\Delta
_{i}^{m}}(u) 1_{\Delta _{j}^{m}}(v) dvdu \\
&=&\sum _{i,j=1; i\not= j} ^{m} \frac{1}{ \vert \Delta _{i}^{m}
\vert \vert \Delta _{j}^{m} \vert } \left( \int _{\Delta _{i}^{m}}
\int _{\Delta _{j}^{m}}F(t,x,y) dydx \right) \\
&&\times \sum _{k,l=1; k\not=l}^{[nt]} n^{2}\frac{\xi
_{k}}{\sqrt{n}} \frac{\xi_{l}}{\sqrt{n}} \lambda \left(
[\frac{k-1}{n}, \frac{k}{n}) \bigcap \Delta _{i}^{m} \right)\lambda
\left( [\frac{l-1}{n}, \frac{l}{n}) \bigcap \Delta _{j}^{m}
\right)\\
&=&\sum _{i,j=1; i\not= j} ^{m} \frac{1}{ \vert \Delta _{i}^{m}
\vert \vert \Delta _{j}^{m} \vert } \left( \int _{\Delta _{i}^{m}}
\int _{\Delta _{j}^{m}}F(t,x,y) dydx \right) \\
&&\times \sum _{k;  [\frac{k-1}{n}, \frac{k}{n}) \subset \Delta
_{i}^{m} }\sum_{l\not= k; [\frac{l-1}{n}, \frac{l}{n}) \subset
\Delta _{j}^{m}} \frac{\xi _{k}}{\sqrt{n}} \frac{\xi_{l}}{\sqrt{n}}
\end{eqnarray*}
and on the other hand  by using (\ref{zpim}) and (\ref{rw}), one has
\begin{equation*}
Z^{\pi ^{m}, n} _{t} =\sum _{i,j=1; i\not= j} ^{m} \frac{1}{ \vert
\Delta _{i}^{m} \vert \vert \Delta _{j}^{m} \vert } \left( \int
_{\Delta _{i}^{m}} \int _{\Delta _{j}^{m}}F(t,x,y) dydx \right) \sum
_{k =[nt_{i-1}^{m}]+1}^{[nt_{i}^{m}]} \sum _{l=
[nt_{j-1}^{m}]+1}^{[nt_{j}^{m}]} \frac{\xi _{k}}{\sqrt{n}}
\frac{\xi_{l}}{\sqrt{n}}
\end{equation*}
and it is not difficult to see that $Z^{\pi ^{m}, n} _{t}$ and
$Z^{',\pi ^{m}, n} _{t}$  coincide.

\vskip0.3cm

{\it Step 4: } At this point we conclude that the family of
processes  $Z^{',n}$ converges in the sense of finite dimensional
distributions to the Rosenblatt process $Z_{t}$. Let $h$ be a
function defined on $\mathbb{R}^{p} $ and consider $s_{1}, \ldots ,
s_{p}\in [0,T]$. We will show that
$$\mathbb{E} \left( h( Z^{',n}_{s_{1}}, \ldots , Z^{',n}_{s_{p}})\right) -\mathbb{E}\left( h(Z_{s_{1}}, \ldots
, Z_{s_{p}} )\right)  $$ converges to zero as $n\to \infty$. This
can be bounded by $A+B+C$ where
$$A= \left| \mathbb{E}\left( h(Z_{s_{1}}, \ldots
, Z_{s_{p}} )\right) -\mathbb{E}\left( h(Z^{\pi ^{m}}_{s_{1}},
\ldots , Z^{\pi ^{m}}_{s_{p}} )\right)\right| $$

$$B= \left| \mathbb{E}\left( h(Z^{\pi ^{m},n }_{s_{1}}, \ldots
, Z^{\pi ^{m},n}_{s_{p}} )\right) -\mathbb{E}\left( h(Z^{\pi
^{m}}_{s_{1}}, \ldots , Z^{\pi ^{m}}_{s_{p}} )\right)\right|$$ and
$$C=  \left| \mathbb{E}\left( h(Z^{\pi ^{m},n }_{s_{1}}, \ldots
, Z^{\pi ^{m},n}_{s_{p}} )\right) -\mathbb{E}\left(
h(Z^{',n}_{s_{1}}, \ldots , Z^{',n}_{s_{p}} )\right)\right|.$$ By
Step 1 and Step 3 we have that for fixed $c>0$ the terms $A$ and $C$
are bounded (uniformly in $n$) by $c$ when $m $ is large enough. the
term $B$ tends to zero as $n \to \infty$ from Step 2.

\vskip0.5cm

{\it Step 5: } Clearly the family $Z^{',n }$ can be replaced by the
family $Z^{n} $(\ref{zn}) because their difference goes to zero in
$L^{2}(\Omega)$ as $n\to \infty. $ \qed

\vskip0.5cm

Next, we show the tightness.
\begin{prop}\label{tight}
The family $Z^{n}$ given by (\ref{zn}) is  tight.
\end{prop}
{\bf Proof: } Let $s<t$, $s,t\in [0,T]$.  It holds, since the kernel
$F(s,u,v) $ vanishes when $u$ or $v$ are bigger that $s$,
\begin{eqnarray*}
&&\mathbb{E}\left| Z^{n}_{t} -Z^{n}_{s} \right| ^{2} \\
&=&\mathbb{E}\left(   \sum _{i,j=1; i\not= j}
^{[nt]}n^{2}\int_{\frac{i-1}{n} }^{\frac{i}{n}}\int_{\frac{j-1}{n}}^
{\frac{j}{n}} \left[ F\left( \frac{[nt]}{n}, u,v\right) -F\left(
\frac{[ns]}{n}, u,v\right) \right] dvdu \frac{\xi
_{i}}{\sqrt{n}}\frac{\xi
_{j}}{\sqrt{n}}\right) ^{2}\\
&=&  \sum _{i,j=1; i\not= j} ^{[nt]}n^{2}\left( \int_{\frac{i-1}{n}
}^{\frac{i}{n}}\int_{\frac{j-1}{n}}^ {\frac{j}{n}}\left[ F\left(
\frac{[nt]}{n}, u,v\right) -F\left( \frac{[ns]}{n}, u,v\right)
\right]  dvdu  \right) ^{2}\\
&\leq &  \sum _{i,j=1; i\not= j} ^{[nt]}\int_{\frac{i-1}{n}
}^{\frac{i}{n}}\int_{\frac{j-1}{n}}^ {\frac{j}{n}}\left[ F\left(
\frac{[nt]}{n}, u,v\right) -F\left( \frac{[ns]}{n}, u,v\right)
\right] ^{2}dvdu\\
&\leq & \int_{0}^{T} \int_{0}^{T}\left[ F\left( \frac{[nt]}{n},
u,v\right) -F\left( \frac{[ns]}{n}, u,v\right)
\right] ^{2} dvdu\\
&=& \left| \frac{[nt]}{n} - \frac{[ns]}{n} \right| ^{2} .
\end{eqnarray*}
Now the conclusion follows by using exactly the same arguments as in
\cite{Sot}, end of the proof of Theorem 1. \qed

\vskip0.5cm

The main result of this section is a consequence of Proposition
\ref{findim} and Proposition \ref{tight}.
\begin{theorem}
\label{main} The family of  stochastic processes (\ref{zn})
converges weakly, in the Skorohod topology, to the Rosenblatt
process.
\end{theorem}

\section{Binary market model based on the Rosenblatt process}

The binary market constitutes a financial model where the asset are
traded at discrete times. In classical cases (for example when the
driven process is the Wiener process) the binary model approximates
the Black and Scholes model.

Let us start by introducing the Black and Scholes model   driven by
the Rosenblatt process.  As usually, we will consider two assets: a
safe investment satisfying
\begin{equation}
\label{asset1} B_{t} = \int_{0} ^{t} r_{s} B_{s} ds,
\end{equation}
where $r$ is a differentiable deterministic function and a risky
asset with price dynamic following the stochastic equation
\begin{equation}
\label{asset2} S_{t}= S_{0}+\int_{0}^{t} a_{s} S_{s}ds   +  \sigma
\int_{0}^{t} S_{s} dZ_{s},
\end{equation}
where $Z$ is a Rosenblatt process, $\sigma >0$ and $a$ is a
differentiable deterministic function. The integral with respect to
$Z$ is understood here in a pathwise  sense. Since the trajectories
of the Rosenblatt process are enough regular (in fact, they are
H\"older continuous of order $\delta <H$ and $H>\frac{1}{2}$) it is
possible to consider  pathwise integrals with respect to it and to
solve some stochastic equations in the pathwise sense. We refer,
among others, to \cite{Nou}, \cite{Za} or \cite{RV1}. In particular,
the solution of (\ref{asset2}) is given by (see \cite{NouTud},
\cite{Za}, \cite{Nou}
\begin{equation}
\label{sol2} S_{t}= S_{0} e^{\int_{0}^{t} a_{s}ds + \sigma Z_{t}
}, \hskip0.5cm t\in [0,T].
\end{equation}
Clearly the solution of (\ref{asset1}) is
\begin{equation}
\label{sol1} B_{t}= B_{0} e^{\int _{0}^{t} r_{s}ds }, \hskip0.5cm
t\in [0,T].
\end{equation}
Moreover, we will assume in the sequel that the interest rates $a$
and $r$ are deterministic bounded functions.

\vskip0.5cm

Let us describe now the binary market model with Rosenblatt
influence. The two assets are traded now at successive times periods
$t_{1}<t_{2} <\ldots < ...$ and their dynamics are given by
$$B_{n} = (1+ r_{n}) B_{n-1} $$
and
$$S_{n}= (a_{n} + (1+ X_{n}) ) S_{n-1};$$
That means that $B_{n} $ and $S_{n}$ represent the prices of the
bond and of the stock in the period between $t_{n}$ and $t_{n+1}$
and $r_{n}$ and $a_{n}$ are the interest rates valuable in this
period. The stochastic process $X$ is binary, that is, given
$X_{n-1}$ it can takes at time $n $ two possible values denoted by
$d_{n}$("down") and $u_{n}$ ("up"). The binary market excludes
arbitrage opportunities if for every $n$ it holds that (see
\cite{Sot})
\begin{equation}
\label{arb} d_{n}<r_{n}-a_{n}< u_{n}.
\end{equation}

In the following we will choose our binary model to be determined by
\begin{equation}\label{xn}
X_{n}= \Delta Z^{N} _{\frac{n}{N}}
\end{equation}
where $Z^{N}$ is defined by (\ref{zn}) and
\begin{equation}\label{ar}
r_{n}= \frac{1}{N} r_{\frac{n}{N} } \mbox{ and } a_{n}= \frac{1}{N}
a_{\frac{n}{N} }
\end{equation}
where $a$ and $r$ are the interest rates appearing in (\ref{asset1})
and (\ref{asset2}).

 \vskip0.5cm

We have \begin{prop} The binary market model with $X$ $a$ and $r$
given by (\ref{xn}) and (\ref{ar}) converges as $N\to \infty$ to the
Black and Scholes model given by (\ref{asset1}) and (\ref{asset2}).
\end{prop}
{\bf Proof: }Let us consider the jump
\begin{equation*}
\Delta Z _{t} ^{N} =Z^{N}_{t}-Z^{N}_{t-}
\end{equation*}
and the quadratic variation

\begin{equation*}
  [Z^{N}] _{t} = \sum _{s\leq t} (\Delta Z^{N} _{s}) ^{2}.
\end{equation*}
We will show that the process $[Z^{N}]$ converges in $L^{1} ([0,T]
\times \Omega)$ to zero. Then the conclusion will follow exactly as
in \cite{Sot}, proof of Lemma 1.

We have, since the jumps are at times $\frac{k}{N}$, $k$ integer,
$$\mathbb{E}\left| \Delta Z_{t} ^{N}\right| ^{2} \leq \mathbb{E}\left| Z_{t}
-Z_{t-\frac{1}{N}}\right| ^{2} \leq \frac{1}{N^{2H}} $$ and then
$$\mathbb{E}[Z^{N}] _{t} \leq Nt\frac{1}{N^{2H}}=tN^{1-2H}.$$
This implies that
$$\int_{0}^{T}[Z^{N}] _{s}ds \leq c(T) N^{1-2H}$$ which goes to $0$ as $N$ goes to $\infty$. \qed

\vskip0.5cm

The next step is  to show that the market admits arbitrage
opportunities. Clearly, we have

\begin{equation*}
X_{n}=\sigma N \sum _{i,j=1; i\not= j}^{n} \left( \int
_{\frac{i-1}{N}}^{\frac{i}{N}}\int
_{\frac{j-1}{N}}^{\frac{j}{N}}\left( F\left( \frac{n}{N},u,v\right)
-F\left( \frac{n-1}{N},u,v\right)\right) dvdu\right) \xi _{i}\xi
_{j}
\end{equation*}
and we will take the random variables $\xi$ to be binary, that is
$$P(\xi =1) = P(\xi _{i} =-1) =\frac{1}{2} $$
for every $i\geq 1$.  We can write, by isolating the part involving
$\xi_{n}$,
\begin{equation*}
X_{n}=f_{n-1} (\xi _{1}, \ldots ,\xi_{n-1}) + \xi _{n} g_{n-1} (\xi
_{1}, \ldots ,\xi_{n-1})
\end{equation*}
where for every $n\geq2$
\begin{equation*}
f_{n-1}(x_{1}, \ldots , x_{n-1})=\sigma N\sum _{i,j=1; i\not=
j}^{n-1} \left( \int _{\frac{i-1}{N}}^{\frac{i}{N}}\int
_{\frac{j-1}{N}}^{\frac{j}{N}}\left( F\left( \frac{n}{N},u,v\right)
-F\left( \frac{n-1}{N},u,v\right)\right) dvdu\right) x_{i}x_{j}
\end{equation*}
and
\begin{equation*}
g_{n-1} (x_{1}, \ldots , x_{n-1}) = 2\sigma N \sum_{i=1}^{n-1}
\left( \int _{\frac{i-1}{N}}^{\frac{i}{N}}\int
_{\frac{n-1}{N}}^{\frac{n}{N}}F\left( \frac{n}{N}, u,v\right)dvdu
\right) x_{i}.
\end{equation*}
Then obviously
\begin{equation}
\label{un} u_{n}=f_{n-1} (\xi _{1}, \ldots ,\xi_{n-1}) + g_{n-1}
(\xi _{1}, \ldots ,\xi_{n-1})
\end{equation}
and
\begin{equation}
\label{dn} d_{n}=f_{n-1} (\xi _{1}, \ldots ,\xi_{n-1}) - g_{n-1}
(\xi _{1}, \ldots ,\xi_{n-1}).
\end{equation}

\vskip0.3cm

 The last result of our paper is the following.

\begin{prop}
The binary market model with (\ref{ar}) and (\ref{xn}) admits
arbitrage.
\end{prop}
{\bf Proof: } Throughout this proof, we will denote by $c(H)$ a
generic constant depending only on $H$. Let us show now that the
condition (\ref{arb}) fails for some $n\geq 2$. We will actually
prove that the sequence
\begin{equation}\label{converg}
f_{n-1}(1, 1, \ldots , 1) - g_{n-1} (1,1, \ldots 1)\to _{n\to \infty
}\infty
\end{equation}
and then clearly (\ref{arb}) does not hold because $r_{n}$ and
$a_{n}$ are assumed to be bounded.

We have
\begin{eqnarray*}
f_{n-1}(1, 1, \ldots , 1)&=&\sigma N\sum _{i,j=1; i\not= j}^{n-1}
\left( \int _{\frac{i-1}{N}}^{\frac{i}{N}}\int
_{\frac{j-1}{N}}^{\frac{j}{N}}\left( F\left( \frac{n}{N},u,v\right)
-F\left( \frac{n-1}{N},u,v\right)\right) dvdu\right) \\
&=&\sigma N\sum _{i,j=1; i}^{n-1} \left( \int
_{\frac{i-1}{N}}^{\frac{i}{N}}\int
_{\frac{j-1}{N}}^{\frac{j}{N}}\left( F\left( \frac{n}{N},u,v\right)
-F\left( \frac{n-1}{N},u,v\right)\right) dvdu\right) \\
&&- \sum _{i=1}^{n-1} \left( \int _{\frac{i-1}{N}}^{\frac{i}{N}}\int
_{\frac{i-1}{N}}^{\frac{i}{N}}\left( F\left( \frac{n}{N},u,v\right)
-F\left( \frac{n-1}{N},u,v\right)\right) dvdu\right)\\
&:=& \sigma N(A-B).
\end{eqnarray*}
Using the expression of the kernel $F$ and (\ref{dK}), the term $A$
can be minorized as follows
\begin{eqnarray*}
A&=&  c(H)\int _{0}^{\frac{n-1}{N}}\int _{0}^{\frac{n-1}{N}}
 \int _{\frac{n-1}{N}} ^{\frac{n}{N}} (a-u)^{H'-\frac{3}{2}}(a-v)^{H'-\frac{3}{2}}
a^{2H' -1} u^{\frac{1}{2}-H'}v^{\frac{1}{2}-H'}dadvdu\\
&\geq & c(H)\left( \frac{n-1}{N} \right) ^{1-2H'}\int
_{0}^{\frac{n-1}{N}}\int _{0}^{\frac{n-1}{N}}
 \int _{\frac{n-1}{N}} ^{\frac{n}{N}} (a-u)^{H'-\frac{3}{2}}(a-v)^{H'-\frac{3}{2}}
a^{2H' -1} dadvdu\\
&\geq & c(H)\left( \frac{n-1}{N} \right) ^{1-2H'}\int
_{\frac{n-1}{N}} ^{\frac{n}{N}}da \left( a^{H'-\frac{1}{2}} -
(a-\frac{n-1}{N} )
^{H'-\frac{1}{2}}\right)^{2} a^{2H'-1}\\
&\geq & c(H)\left( \frac{n-1}{N} \right) ^{1-2H'}\int
_{\frac{n-1}{N}}
^{\frac{n}{N}}da  \left(\frac{ (n-1)^{H'-\frac{1}{2}}-1}{N^{H'-\frac{1}{2}}} \right) ^{2} a^{2H'-1}\\
&\geq &c(H) \frac{1}{N^{2H'}}
(n-1)^{1-2H'}((n-1)^{H'-\frac{1}{2}}-1)^{2}(n^{2H'}-(n-1)^{2H'} )\\
&=&c(H)\frac{1}{N^{2H'}}o(n^{2H' -1}).
\end{eqnarray*}
We majorize now the term $B$. We can write
\begin{eqnarray*}
B&\leq & c(H)\sum _{i=1}^{n-1}  \int
_{\frac{i-1}{N}}^{\frac{i}{N}}\int
_{\frac{i-1}{N}}^{\frac{i}{N}}\int _{\frac{n-1}{N}} ^{\frac{n}{N}}
(a-u)^{H'-\frac{3}{2}}(a-v)^{H'-\frac{3}{2}}
a^{2H' -1} u^{\frac{1}{2}-H'}v^{\frac{1}{2}-H'}dadvdu\\
&\leq& c(H) \left( \frac{n}{N}\right) ^{2H' -1} \sum
_{i=1}^{n-1}\int _{\frac{i-1}{N}}^{\frac{i}{N}}\int
_{\frac{i-1}{N}}^{\frac{i}{N}}\int _{\frac{n-1}{N}} ^{\frac{n}{N}} (a- \frac{i}{N})^{2H' -3} u^{\frac{1}{2}-H'}v^{\frac{1}{2}-H'}dadvdu\\
&\leq & c(H)\left( \frac{n}{N}\right) ^{2H' -1} \sum
_{i=1}^{n-1}\int _{\frac{i-1}{N}}^{\frac{i}{N}}\int
_{\frac{i-1}{N}}^{\frac{i}{N}}u^{\frac{1}{2}-H'}v^{\frac{1}{2}-H'}\frac{
(n-i) ^{2H'-2} -(n-(i-1) ) ^{2H'-2}}{N^{2H'-2} }dvdu \\
&\leq & c(H)\left( \frac{n}{N}\right) ^{2H' -1}\sum _{i=1}^{n-1}
\left(  \left( \frac{i}{N} \right) ^{\frac{3}{2}-H'}-\left(
\frac{i-1}{N} \right) ^{\frac{3}{2}-H'}\right) ^{2}\frac{ (n-i)
^{2H'-2} -(n-(i-1) ) ^{2H'-2}}{N^{2H'-2} }\\
&\leq & c(H)\left( \frac{n}{N}\right) ^{2H' -1} \frac{1}{N^{3-2H'}}
\left( \frac{n-1}{N}\right) ^{2H'-2}\\
&\leq& c(H) \frac{1}{N^{2H'}} o(n^{4H'-3}).
\end{eqnarray*}
From the above computations we obtain that $f_{n-1}(1, \ldots ,1)$
converges to $\infty $ as $n\to \infty$ because $2H'-1 >4H'-3$. Let
us estimate the term $g_{n-1}(1, \ldots , 1)$.

\begin{eqnarray*}
&& g_{n-1}(1, \ldots , 1)\\
&\leq&  2\sigma N \sum_{i=1}^{n-1} \left( \int
_{\frac{i-1}{N}}^{\frac{i}{N}}\int
_{\frac{n-1}{N}}^{\frac{n}{N}}F\left( \frac{n}{N}, u,v\right)dvdu
\right)\\
&\leq & 2\sigma N c(H)\left( \frac{n}{N}\right) ^{2H'
-1}\sum_{i=1}^{n-1}  \int _{\frac{i-1}{N}}^{\frac{i}{N}}\int
_{\frac{n-1}{N}}^{\frac{n}{N}}\int _{v}^{\frac{n}{N}}
u^{\frac{1}{2}-H'}v^{\frac{1}{2}-H'}(\frac{n-1}{N}-u)
^{H'-\frac{3}{2}} (a-v) ^{H'-\frac{3}{2}}dadvdu \\
&\leq&2\sigma N c(H)\left( \frac{n}{N}\right) ^{2H'
-1}\left(\frac{n-1}{N} \right) ^{\frac{1}{2}-H'}\sum_{i=1}^{n-1}\int
_{\frac{i-1}{N}}^{\frac{i}{N}}\int
_{\frac{n-1}{N}}^{\frac{n}{N}}u^{\frac{1}{2}-H'}(\frac{n-1}{N}-u) ^
{H'-\frac{3}{2}} (\frac{n}{N} -v) ^{H'-\frac{1}{2}}dvdu \\
&\leq &2\sigma N c(H)\left( \frac{n}{N}\right) ^{2H'
-1}\left(\frac{n-1}{N} \right)^{\frac{1}{2}-H'}\int
_{0}^{\frac{n-1}{N}}u^{\frac{1}{2}-H'}(\frac{n-1}{N}-u) ^
{H'-\frac{3}{2}} du\\
&\leq &2\sigma N c(H) \frac{1}{N^{2H' -1}} o(n^{H'-\frac{1}{2}} )
\end{eqnarray*}
because the integral $\int
_{0}^{\frac{n-1}{N}}u^{\frac{1}{2}-H'}(\frac{n-1}{N}-u) ^
{H'-\frac{3}{2}} du$ is equal to $\beta (\frac{3}{2}-H',
H'-\frac{1}{2})$. We then obtain (\ref{converg}). \qed

\vskip0.5cm

{\bf Comments: } i) A concrete arbitrage opportunity can be easily
described. For example, suppose that $a>r$. Ar a certain time
$n_{0}$ we have $d_{n_{0}}>0$ because of (\ref{converg}); suppose
that the stock price was increasing up to this time $n_{0}$. Then,
buy $M$ stocks and your wealth at time $n_{0}+1$ will be positive
since $MS_{n_{0}+1}>MS_{n_{0}}$.

\vskip0.3cm

 ii) We cannot expect to have no-arbitrage when $H\in
\frac{1}{2}$ as in the fractional Brownian motion case because now
the limit process at $H=\frac{1}{2}$ is not necessarily a
martingale.

\newpage

\section{Appendix: Monte Carlo Simulation}\label{Ap}

This simulation study is intended to show our proposed model using
simulated data. We simulate data using different values for $ H$.
Below, we describe the procedure used in generating the data to be
used in the simulation study.

\begin{figure}[!htp]
\begin{minipage}[b]{.3\textwidth}
\includegraphics[scale=.35]{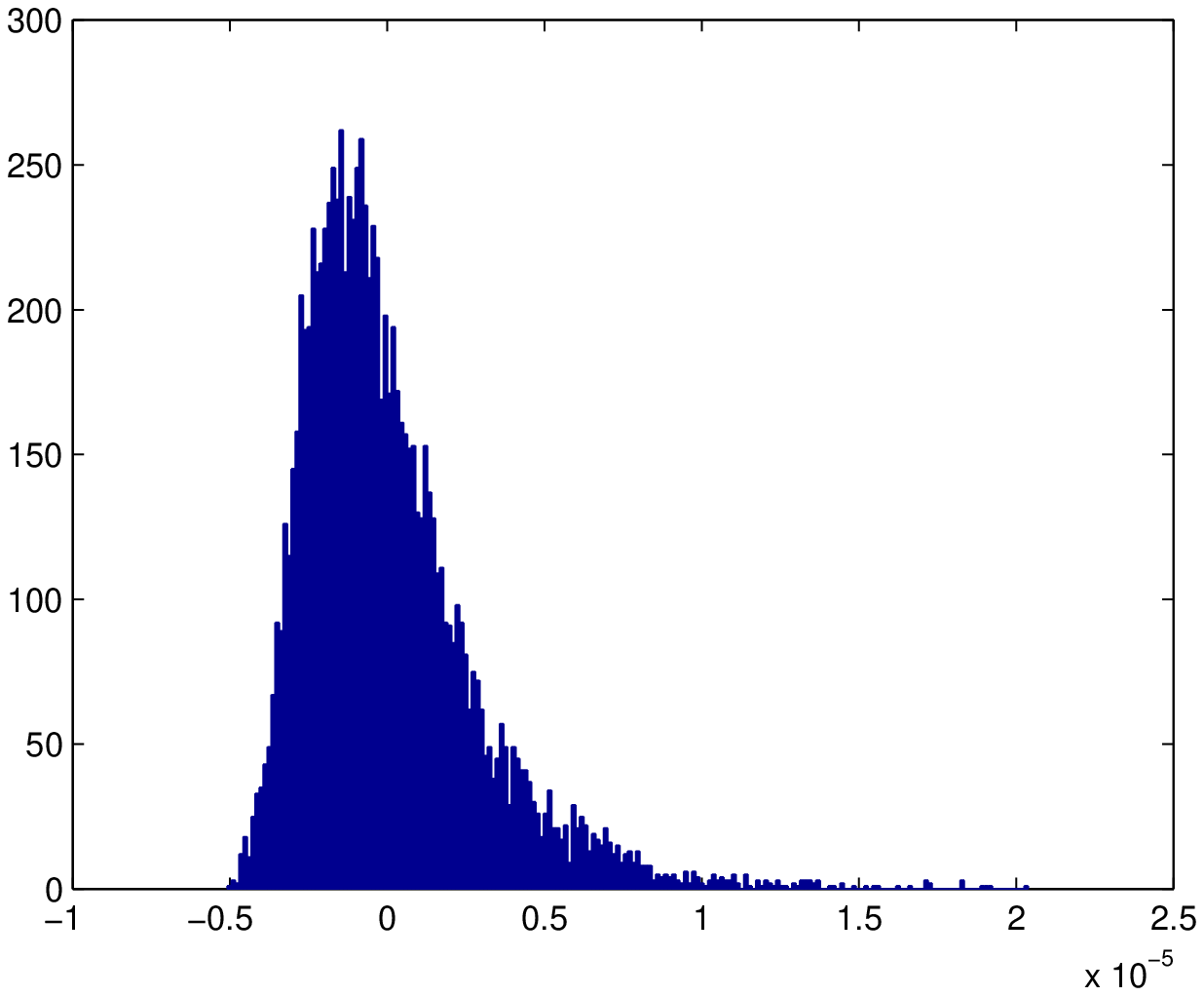}
\end{minipage}
\begin{minipage}[b]{.3\textwidth}
\includegraphics[scale=.35]{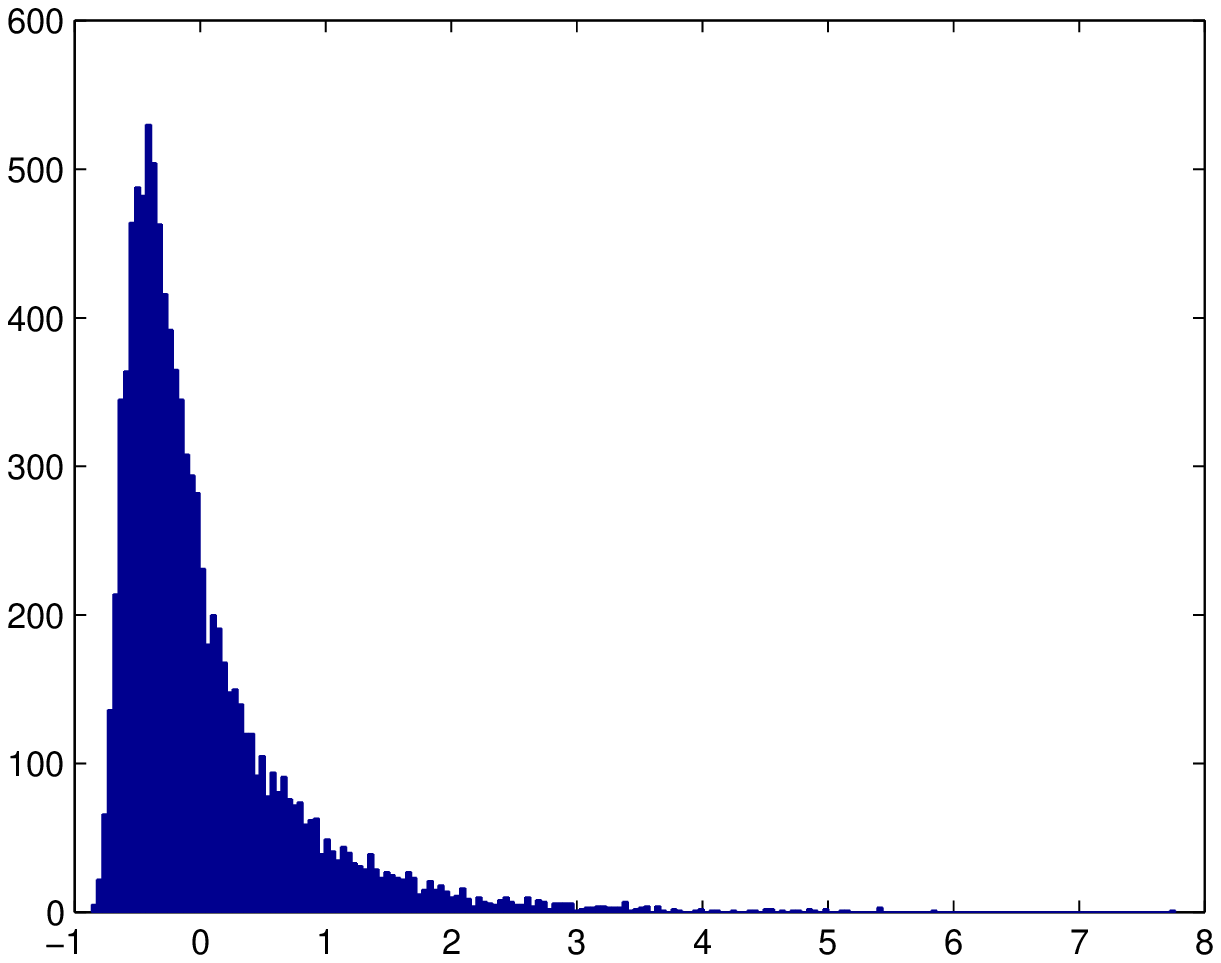}
\end{minipage}
\begin{minipage}[b]{.3\textwidth}
\includegraphics[scale=.35]{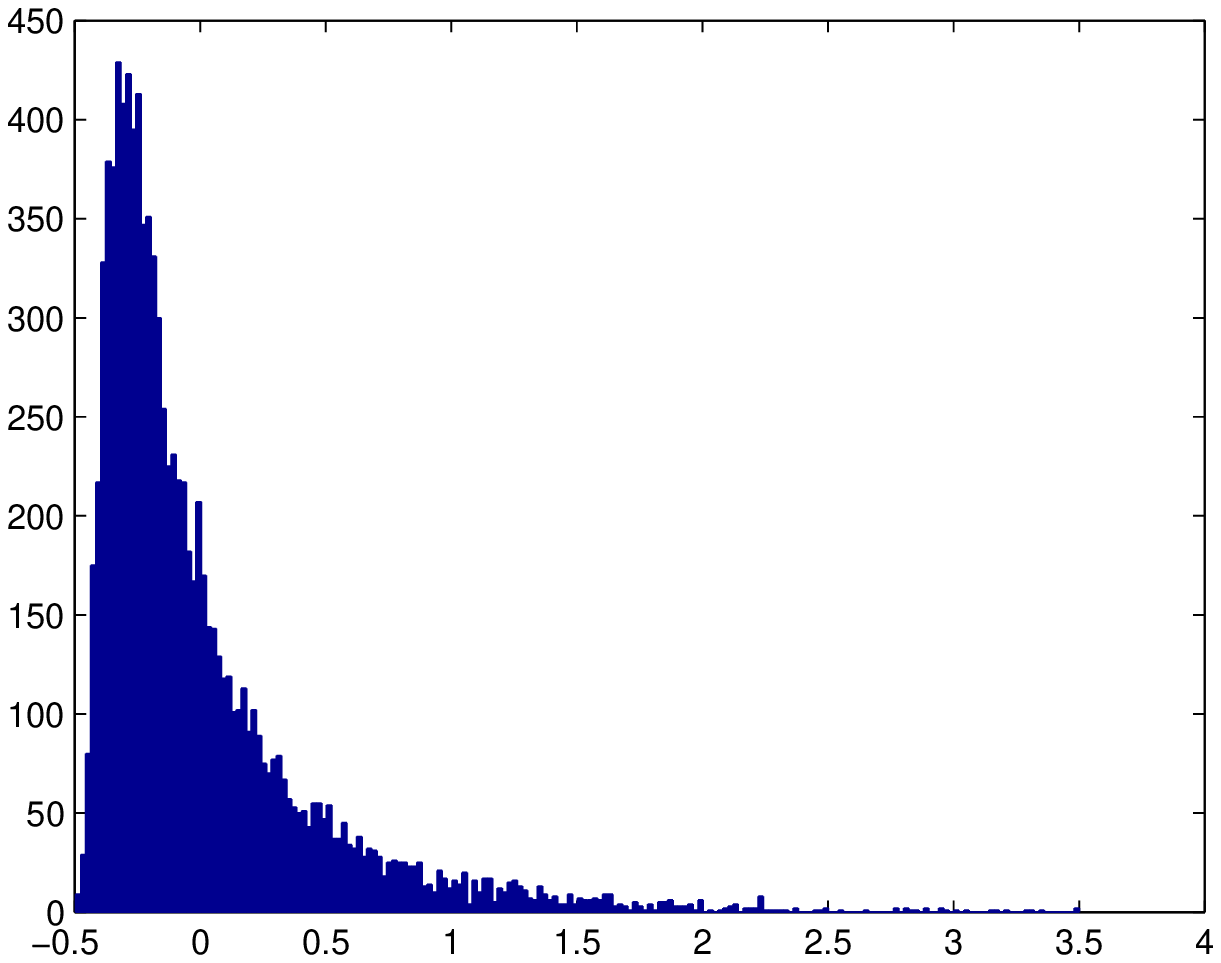}
\end{minipage}
\caption{Marginal distributions for $H=0.5$; $H=0.8$ and $H=0.9$.}
\label{fig1:hist}
\end{figure}

We have implemented this simulation on a standard personal
computing platform (PC), and have observed that it performs very
well using simulated data as can be seen from the simulated data
in the figures 1 and 2 below. Despite the apparent algebraic
complexity of the equations (\ref{zn}), the problem poses no
difficulty for standard symbolic algebra packages. Using Matlab's
simulations and algebra capabilities yielded the best computing
times. In our implementation, which performs an iteration of the
algorithm from $i=0$ to $i=n$. Figure \ref{fig1:hist} shows the
histogram for a fixed time $t$ for the marginal density of the
Rosenblatt process. We can see the skew structure of the
distribution.

\begin{figure}[!htp]
\begin{minipage}[b]{.5\textwidth}
\includegraphics[scale=.45]{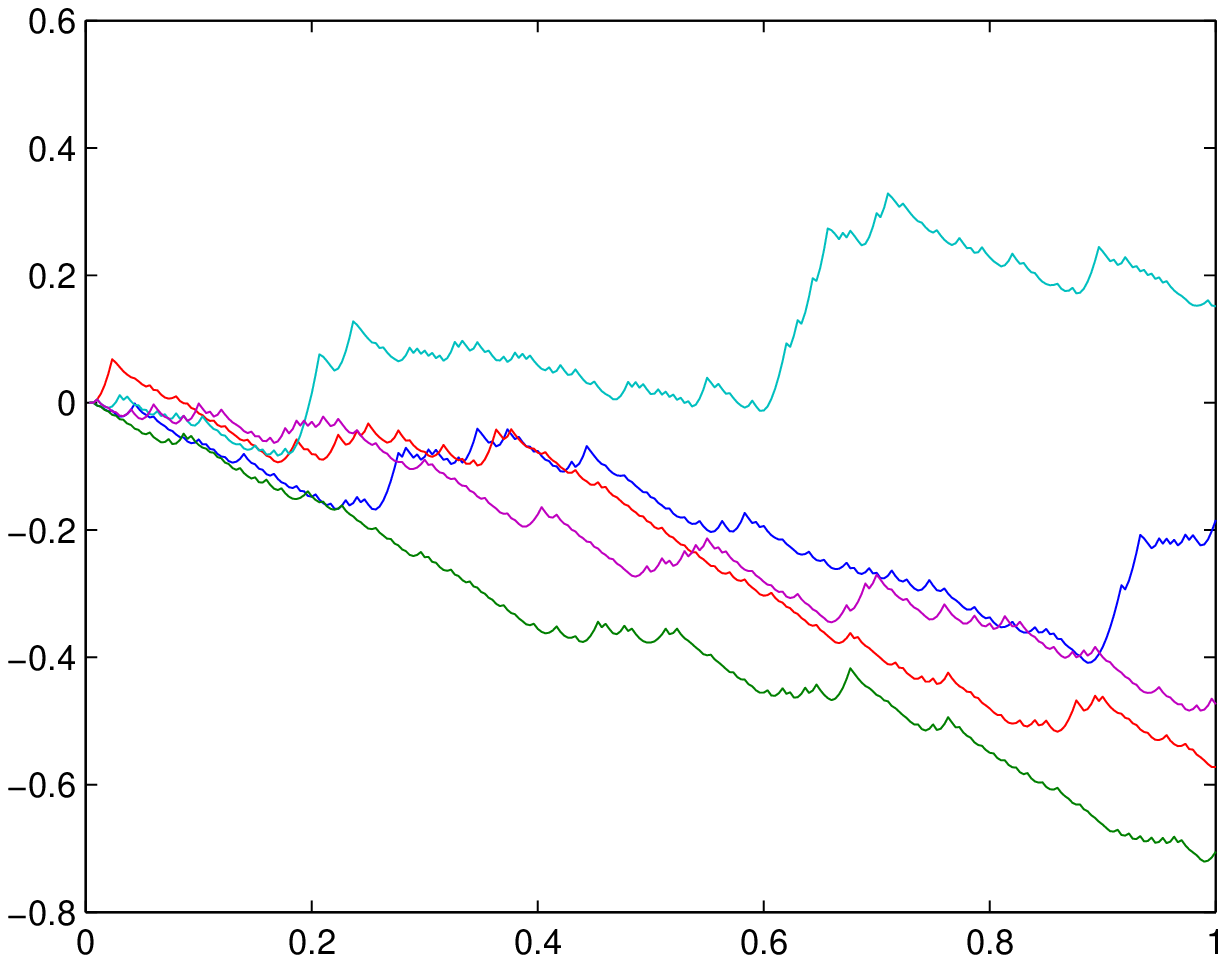}
\end{minipage}
\begin{minipage}[b]{.3\textwidth}
\includegraphics[scale=.45]{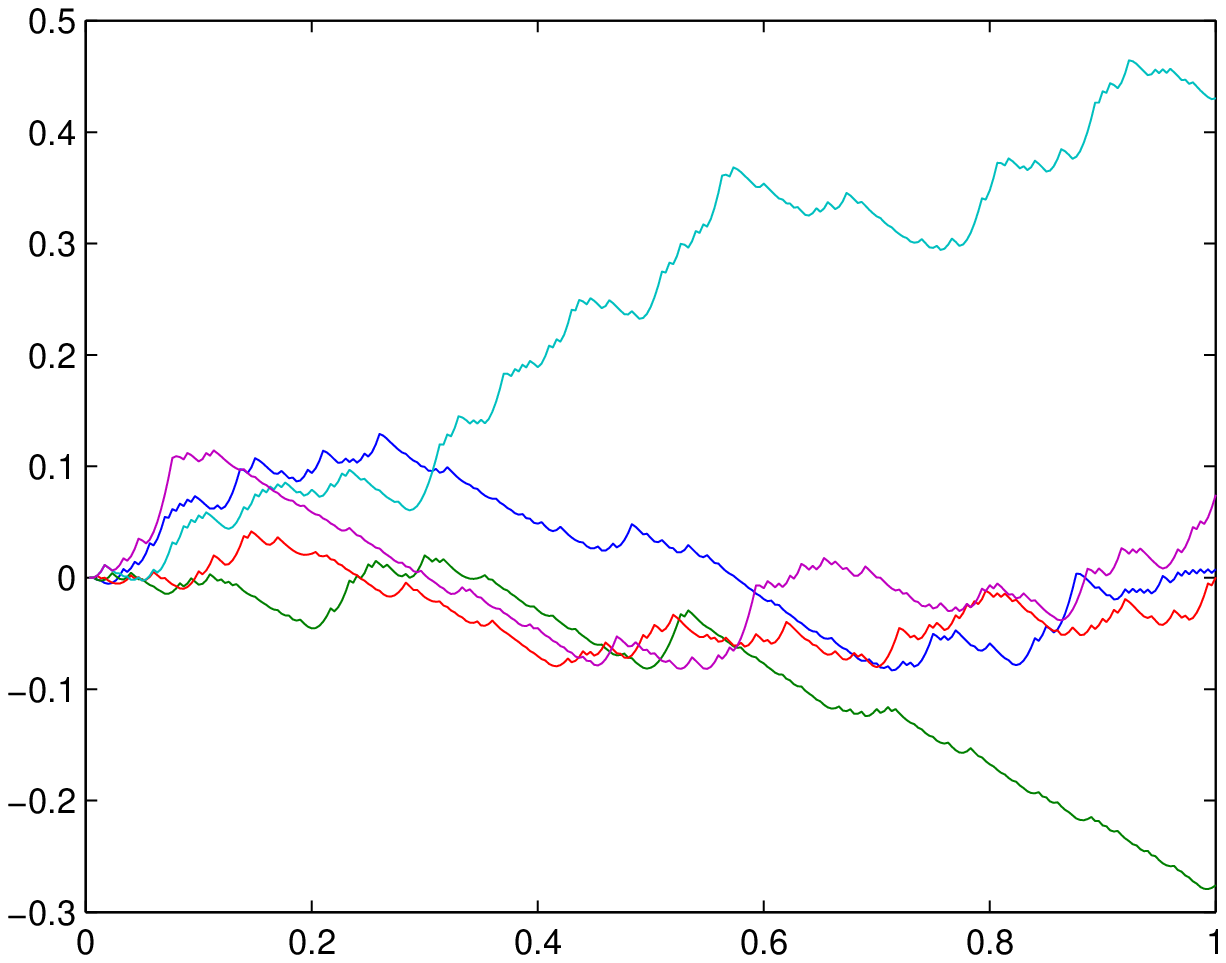}
\end{minipage}
\caption{Simulations of the Trajectory for the Rosenblatt process
with $H=0.8$ and $H=0.9$.} \label{fig2:paths}
\end{figure}

Figure \ref{fig2:paths} shows the some paths of the discretization
for the Rosenblatt process. We use the values for the parameter
$H$ ($H=0.8$ and $H=0.9$).

{\bf Acknowledgments:} This work was partially supported by the
research project Nucleus Millenium P04-069-F "Information and
Randomness: Fundamentals and Applications; Laboratories in
Mathematics of Genome and Stochastic Simulation". The first author
was supported partially by the research project Fondecyt Reg. N.
1050843, Chile and Proyecto Anillo ACT-13: "Laboratorio ANESTOC".

\end{document}